\documentclass[12pt]{amsart}

\usepackage[hypertexnames=false, colorlinks, citecolor=red, linkcolor=red]{hyperref} 
\usepackage{amstext}
\usepackage{amsthm}
\usepackage{amsrefs}
\usepackage{amsmath}
\usepackage{amssymb}
\usepackage{latexsym}
\usepackage{amsfonts}

\input txdtools

\begin{document}

\newcommand{\ci}[1]{_{ {}_{\scriptstyle #1}}}

\newcommand{\norm}[1]{\ensuremath{\|#1\|}}
\newcommand{\abs}[1]{\ensuremath{\vert#1\vert}}
\newcommand{\p}{\ensuremath{\partial}}
\newcommand{\pr}{\mathcal{P}}

\newcommand{\pbar}{\ensuremath{\bar{\partial}}}
\newcommand{\db}{\overline\partial}
\newcommand{\D}{\mathbb{D}}
\newcommand{\B}{\mathbb{B}}
\newcommand{\Sp}{\mathbb{S}}
\newcommand{\T}{\mathbb{T}}
\newcommand{\R}{\mathbb{R}}
\newcommand{\Z}{\mathbb{Z}}
\newcommand{\C}{\mathbb{C}}
\newcommand{\N}{\mathbb{N}}
\newcommand{\scrH}{\mathcal{H}}
\newcommand{\scrL}{\mathcal{L}}
\newcommand{\td}{\widetilde\Delta}

\newcommand{\La}{\langle }
\newcommand{\Ra}{\rangle }
\newcommand{\rk}{\operatorname{rk}}
\newcommand{\card}{\operatorname{card}}
\newcommand{\ran}{\operatorname{Ran}}
\newcommand{\osc}{\operatorname{OSC}}
\newcommand{\im}{\operatorname{Im}}
\newcommand{\re}{\operatorname{Re}}
\newcommand{\tr}{\operatorname{tr}}
\newcommand{\vf}{\varphi}
\newcommand{\f}[2]{\ensuremath{\frac{#1}{#2}}}


\newcommand{\entrylabel}[1]{\mbox{#1}\hfill}

\newenvironment{entry}
{\begin{list}{X}%
  {\renewcommand{\makelabel}{\entrylabel}%
      \setlength{\labelwidth}{55pt}%
      \setlength{\leftmargin}{\labelwidth}
      \addtolength{\leftmargin}{\labelsep}%
   }%
}%
{\end{list}}


\numberwithin{equation}{section}

\newtheorem{thm}{Theorem}[section]
\newtheorem{lm}[thm]{Lemma}
\newtheorem{cor}[thm]{Corollary}
\newtheorem{conj}[thm]{Conjecture}
\newtheorem{prob}[thm]{Problem}
\newtheorem{prop}[thm]{Proposition}
\newtheorem*{prop*}{Proposition}

\theoremstyle{remark}
\newtheorem{rem}[thm]{Remark}
\newtheorem*{rem*}{Remark}

\title{Multi-Parameter Div-Curl Lemmas}

\author[M. T. Lacey]{Michael T. Lacey$^{\ddagger}$}
\address{Michael T. Lacey, School of Mathematics\\ Georgian Institute of Technology\\ 686 Cherry Street\\ Atlanta, GA USA 30332--0160}
\email{lacey@math.gatech.edu}
\thanks{$\ddagger$ Research supported in part by a National Science Foundation DMS grant.}

\author[S. Petermichl]{Stefanie Petermichl}
\address{Stefanie Petermichl, UniversitŽ Paul Sabatier\\ Institut de MathŽmatiques de Toulouse\\
118 route de Narbonne\\ F-31062 Toulouse Cedex 9, France }
\email{stefanie.petermichl@math.univ-toulouse.fr}

\author[J. C. Pipher]{Jill C. Pipher$^{\dagger}$}
\address{Jill C. Pipher, Department of Mathematics\\ Brown University \\ 151 Thayer Street\\ Providence, RI USA 02912}
\email{pipher@math.brown.edu}
\thanks{$\dagger$ Research supported in part by a National Science Foundation DMS grant.}

\author[B. D. Wick]{Brett D. Wick$^{\ast}$}
\address{Brett D. Wick, School of Mathematics\\ Georgian Institute of Technology\\ 686 Cherry Street\\ Atlanta, GA USA 30332--0160}
\address{\hspace*{1.1in} Institute for Mathematics\\ University of Paderborn\\ Warburger Str. 100\\ 33098 Paderborn Germany}
\email{wick@math.gatech.edu}
\urladdr{http://people.math.gatech.edu/$\sim$bwick6/}
\thanks{$\ast$ Research supported in part by a National Science Foundation DMS grant and an Alexander von Humboldt Research Fellowship.}

\subjclass[2000]{Primary: }

\begin{abstract} 
We study the possible analogous of the Div-Curl Lemma in classical harmonic analysis and partial differential equations, but from the point of view of the multi-parameter setting.  In this context we see two possible Div-Curl lemmas that arise.  Extensions to differential forms are also given.   
\end{abstract}

\maketitle

\section{Introduction}
Suppose that $E, B\in L^2(\R^n;\R^n)$ are vector-fields on $\R^n$ taking values in $\R^n$.  Simply by supposing that the vector-fields belong to $L^2$, one immediately sees that their dot product $E\cdot B\in L^1(\R^n)$ with $\norm{E\cdot B}_{L^1(\R^n)}\leq\norm{E}_{L^2(\R^n;\R^n)}\norm{B}_{L^2(\R^n;\R^n)}$.  If in addition the vector fields satisfy the additional conditions that
\begin{equation*}
\operatorname{div} E(x)=0\quad \operatorname{curl} B(x)=0
\end{equation*}

Then we have that $E\cdot B\in H^1(\R^n)$ with $\norm{E\cdot B}_{H^1(\R^n)}\lesssim\norm{E}_{L^2(\R^n;\R^n)}\norm{B}_{L^2(\R^n;\R^n)}$.  Thus the conditions of vector fields being curl- and  divergence-free induce enough cancellation that $E\cdot B\in H^1(\R^n)$.  Since the proof of this fact motivates what we are interested in, we give the main ideas now, which can be found in the paper by Coifman, Lions, Meyer and Semmes, \cite{CLMS}.

This fact is easy to see by using the commutator theorems of Coifman, Rochberg and Weiss, see \cite{CRW}.  Since $B$ is a curl-free vector field there exists a function $\varphi\in L^2(\R^n)$ such that the components of $B$, $B_j$ are given by $R_j\varphi$ and $\norm{B}_{L^2(\R^n;\R^n)}\approx\norm{\varphi}_{L^2(\R^n)}$.  With this observation we have that
\begin{eqnarray*}
E(x)\cdot B (x) & = & \sum_{j=1}^n E_j(x)B_j(x)\\
& = & \sum_{j=1}^n E_j(x) R_j\varphi(x)+\varphi(x)R_j E_j(x)-\varphi(x)R_jE_j(x)\\
& = & \sum_{j=1}^n E_j(x)R_j\varphi(x)+\varphi(x)R_jE_j(x),
\end{eqnarray*}
with the last equality following because the vector field $E$ is divergence-free and so $\sum_{j=1}^nR_jE_j(x)=0$ (easily seen by taking the Fourier Transform).  Testing this equality over all functions $b\in BMO(\R^n)$ we see that 
\begin{eqnarray*}
\int_{\R^n}E(x)\cdot B (x) b(x)dx & = & \sum_{j=1}^n \int_{\R^n}b(x)\left(E_j(x)R_j\varphi(x)+\varphi(x)R_jE_j(x)\right)dx\\
& = & \sum_{j=1}^n \int_{\R^n}\left[b, R_j\right](E_j)(x)\varphi(x) dx.\\
\end{eqnarray*}
Since $b\in BMO(\R^n)$ we have that each of the above commutators is a bounded operator on $L^2(\R^n)$ with norm controlled by the norm of $b\in BMO(\R^n)$.  Putting this together, we see that
\begin{eqnarray*}
\left\vert\int_{\R^n}E(x)\cdot B (x) b(x)dx\right\vert & \lesssim & \norm{b}_{BMO}\norm{E}_{L^2(\R^n;\R^n)}\norm{\varphi}_{L^2(\R^n)}\approx \norm{b}_{BMO}\norm{E}_{L^2(\R^n;\R^n)}\norm{B}_{L^2(\R^n;\R^n)}.
\end{eqnarray*}
Fefferman duality between $H^1(\R^n)$ and $BMO(\R^n)$ then yields that 
$$
\norm{E\cdot B}_{H^1(\R^n)}\lesssim\norm{E}_{L^2(\R^n;\R^n)}\norm{B}_{L^2(\R^n;\R^n)}.
$$

In the rest of this note we want to address the analogous question in the multi-parameter setting.  This leads to two new results, one inspired by a direct generalization of these ideas, while the other arises by first considering the product structure, and then writing down the conditions that arise.

For the statement of the first theorem, let $\mathcal{M}_{nm}$ denote the set of all $n\times m$ matrices.  Given two matrices $A,B\in\mathcal{M}_{nm}$ we define the ``dot product'' between $A$ and $B$ by 
$$
A\cdot B:=\sum_{j=1,k=1}^{n,m}A_{jk}B_{jk}.
$$
We point out that this is the Hilbert-Schmidt inner product for two matrices and more generally this is referred to as the Schur product of two matrices.

Thus, given $F\in L^2(\R^n\times\R^m;\mathcal{M}_{nm})$, for fixed $j$, $(1\leq j\leq n)$ and $x\in\R^m$, we can view $F_{j,k}(x,y)$ $(1\leq k\leq m)$ as a vector field on $\R^m$ taking values in $\R^m$.  A similar statement applies with the roles of the variables interchanged.  This leads to the following Theorem.

\begin{thm}
\label{product}
Let $1<p<\infty$ and $\frac{1}{p}+\frac{1}{q}=1$.  Suppose that $E\in L^p(\R^n\times \R^m;\mathcal{M}_{nm})$ and $B\in L^q(\R^n\times \R^m;\mathcal{M}_{nm})$ are such that 
\begin{equation*}
\operatorname{div}_x E_{j,k}(x,y)=0\quad \operatorname{curl}_x B_{j,k}(x,y)=0\quad\forall y\in\R^m,\forall k
\end{equation*}
and
\begin{equation*}
\operatorname{div}_y E_{j,k}(x,y)=0\quad \operatorname{curl}_y B_{j,k}(x,y)=0\quad\forall x\in\R^n,\forall j.
\end{equation*}
Then $E\cdot B$ belongs to product $H^1(\R^n\otimes\R^m)$ with
$$
\norm{E\cdot B}_{H^1(\R^n\otimes\R^m)}\lesssim\norm{E}_{L^p(\R^n\times\R^m;\mathcal{M}_{nm})}\norm{B}_{L^q(\R^n\times\R^m;\mathcal{M}_{nm})}.
$$
\end{thm}

The second theorem is the most obvious generalization of the Div-Curl Lemma from one-parameter harmonic analysis to the multi-parameter setting.  This generalization requires that we be working in the ``square'' case since we need both $x,y\in\R^n$.

\begin{thm}
\label{uniformly}
Let $1<p<\infty$ and $\frac{1}{p}+\frac{1}{q}=1$.  Suppose that $E\in L^p(\R^n\times \R^n;\R^n)$ and $B\in L^q(\R^n\times \R^n;\R^n)$ are such that 
\begin{equation*}
\operatorname{div}_x E(x,y)=0\quad \operatorname{curl}_x B(x,y)=0\quad\forall y\in\R^n
\end{equation*}
and
\begin{equation*}
\operatorname{div}_y E(x,y)=0\quad \operatorname{curl}_y B(x,y)=0\quad\forall x\in\R^n.
\end{equation*}
Then
$$
\int_{\R^n}\norm{E(x,\cdot)\cdot B(x,\cdot)}_{H^1(\R^n)}dx\lesssim\norm{E}_{L^p(\R^n\times\R^n;\R^n)}\norm{B}_{L^q(\R^n\times\R^n;\R^n)}
$$
and
$$
\int_{\R^n}\norm{E(\cdot,y)\cdot B(\cdot,y)}_{H^1(\R^n)}dy\lesssim\norm{E}_{L^p(\R^n\times\R^n;\R^n)}\norm{B}_{L^q(\R^n\times\R^n;\R^n)}.
$$
\end{thm}

At this point, we remark that the theorems have extensions to the case of more than two parameters.  For the sake of notation and convenience, we state the theorems and give the proofs only in the case of two parameters.  For the correct generalization of Theorem \ref{product} one must use higher tensors, but the ideas in these cases are identical.  The interested reader can easily extend these results.

\section{Multi-Parameter Div-Curl Lemmas}

In this section we look at the two possible generalizations of the Div-Curl lemma to the product setting.  The first will be obtained by considering a product of divergence- and curl-free vector fields.  In this generalization we end up considering higher tensors that arise because of possible multiplications that arise.

The second generalization we consider vector fields $E$ and $B$ that are uniformly curl- and divergence-free in each variable separately.  With these hypotheses we will be able to show that the resulting function is in fact uniformly $H^1$ in each variable separately.
   
\subsection{Higher Tensor Generalization}  

In this subsection we prove Theorem \ref{product}.  The idea of the proof is to simply mimic the classic one-parameter proof, but to take advantage of the conditions we have on the vector fields $E$ and $B$.

Since $B$ is a curl-free vector field the there exists a function $\varphi\in L^q(\R^n\times \R^m)$ such that the components of $B$, $B_{j,k}$ are given by $R_j^x R_k^y\varphi$ and $\norm{B}_{L^q(\R^n\times\R^m;\mathcal{M}_{nm})}\approx\norm{\varphi}_{L^q(\R^n\times\R^m)}$.  Here we are letting $R_j^x$ denote the $j^{th}$ Riesz transform in the $x$ variable and a similar statement applying to $R_k^y$.  With this observation we have that
\begin{eqnarray*}
E(x,y)\cdot B (x,y) & = & \sum_{j=1}^{n}\sum_{k=1}^{m} E_{j,k}(x,y)B_{j,k}(x,y)\\
& = & \sum_{j=1}^n \sum_{k=1}^m \left\{E_{j,k}(x,y) R_j^x R_k^y\varphi(x,y)\right.\\
 & & +\varphi(x,y)R_j^xR_k^y E_{j,k}(x,y)-\varphi(x,y)R_j^xR_k^yE_{j,k}(x,y)\\
&  & +R_j^xE_{j,k}(x,y)R_k^y\varphi(x,y)-R_j^xE_{j,k}(x,y)R_k^y\varphi(x,y)\\
& &\left. +R_k^yE_{j,k}(x,y)R_j^x\varphi(x,y)-R_k^yE_{j,k}(x,y)R_j^x\varphi(x,y)\right\}\\
& = & \sum_{j=1}^n \sum_{k=1}^m\left\{E_{j,k}(x,y)R_j^xR_k^y\varphi(x,y)+\varphi(x,y)R_j^xR_k^yE_{j,k}(x,y)\right.\\ 
& & \left.+R_j^x\varphi(x,y)R_k^yE_{j,k}(x,y)+R_k^y\varphi(x,y)R_j^xE_{j,k}(x,y)\right\},
\end{eqnarray*}
with the last equality following because the function $E$ is divergence-free in each row and column separately and so $\sum_{j=1}^nR_j^xE_{j,k}(x,y)=0$ for all $k$ and $y\in\R^m$ (easily seen by taking the Fourier Transform).  Similarly, we know that  $\sum_{k=1}^nR_k^yE_{j,k}(x,y)=0$ for all $j$ and $x\in\R^n$. Indeed, we see that the three terms with minus signs from above can be written as (ignoring the sign),
\begin{eqnarray*}
\varphi(x,y)\sum_{j=1}^n\sum_{k=1}^mR_j^xR_k^yE_{j,k}(x,y) & = & \varphi(x,y)\sum_{j=1}^nR_j^x\sum_{k=1}^mR_k^yE_{j,k}(x,y)=0\\
\sum_{j=1}^n\sum_{k=1}^mR_j^xE_{j,k}(x,y)R_k^y\varphi(x,y) & = & \sum_{k=1}^m\left( R_k^y\varphi(x,y)\sum_{j=1}^n R_j^x E_{j,k}(x,y)\right)=0\\
\sum_{j=1}^n\sum_{k=1}^mR_k^yE_{j,k}(x,y)R_j^x\varphi(x,y) & = & \sum_{j=1}^n\left( R_j^x\varphi(x,y)\sum_{k=1}^m R_k^y E_{j,k}(x,y)\right)=0.\\
\end{eqnarray*}
The terms in the parentheses are zero by the hypotheses on the field $E(x,y)$.  Also, the first term is easy to see that it vanishes from the conditions on $E(x,y)$.

Testing this equality over all functions in product BMO, $b\in BMO(\R^n\otimes \R^m)$, and un-ravelling the expression with the Riesz transforms we see that 
\begin{eqnarray*}
\int_{\R^n\times\R^m}E(x,y)\cdot B (x,y) b(x,y)dxdy & = & \sum_{j=1}^n \sum_{k=1}^m\int_{\R^n\times\R^m}\left[\left[b, R_j^x\right],R_k^y\right](E_{j,k})(x,y)\varphi(x,y) dxdy.
\end{eqnarray*}
Since $b\in BMO(\R^n\otimes \R^m)$ we have that each of the above commutators is a bounded operator on $L^p(\R^n\times\R^m)$ with norm controlled by the norm of $b\in BMO(\R^n\otimes\R^m)$, see \cite{LPPW}.  Putting this together, we see that
\begin{eqnarray*}
\left\vert\int_{\R^n\times \R^m}E(x,y)\cdot B (x,y) b(x,y)dxdy\right\vert & \lesssim & \norm{b}_{BMO}\norm{E}_{L^p(\R^n\times\R^m;\mathcal{M}_{nm})}\norm{\varphi}_{L^q(\R^n\times\R^m)}\\
& \approx & \norm{b}_{BMO}\norm{E}_{L^p(\R^n\times\R^m;\mathcal{M}_{nm})}\norm{B}_{L^q(\R^n\times\R^m;\mathcal{M}_{nm})}.
\end{eqnarray*}
Chang-Fefferman duality between product $H^1(\R^n\otimes\R^m)$ and product $BMO(\R^n\otimes\R^m)$, see \cite{CF}, then yields that $\norm{E\cdot B}_{H^1(\R^n\otimes\R^m)}\lesssim\norm{E}_{L^p(\R^n\times\R^m;\mathcal{M}_{nm})}\norm{B}_{L^q(\R^n\times\R^m;\mathcal{M}_{nm})}$.

\subsection{Vector Generalization}
We now turn to the proof of Theorem \ref{uniformly}.  Given a vector-field $F$, we will let $F_j$ denote the $j^{th}$ component of $F$.  Suppose that we have vector fields $E, B\in L^2(\R^n\times \R^n;\R^n)$ with  
\begin{equation*}
\operatorname{div}_x E(x,y)=0\quad \operatorname{curl}_x B(x,y)=0\quad\forall y\in\R^n
\end{equation*}
and
\begin{equation*}
\operatorname{div}_y E(x,y)=0\quad \operatorname{curl}_y B(x,y)=0\quad\forall x\in\R^n.
\end{equation*}

Before, we proof Theorem \ref{uniformly}, we return to the one-parameter case and provided a different method to see the Div-Curl Lemmas in that context.  This will serve as motivation for the multi-parameter setting.

It is easy to see that for a vector-field $F\in L^2(\R^n;\R^n)$ that 
$$
\mathcal{P}_x F(x):=\left(F_1+R_1\left(\sum_{k=1}^n R_kF_k\right),\ldots, F_n+R_n\left(\sum_{k=1}^n R_kF_k\right)\right)
$$ 
projects the vector-field $F$ onto the divergence-free vector fields on $\R^n$.  Indeed, we have that
\begin{eqnarray*}
\sum_{j=1}^n R_j\left(F_j+R_j\left(\sum_{k=1}^nR_kF_k\right)\right)=\sum_{j=1}^n R_jF_j+\sum_{j=1}^n R_j^2\left(\sum_{k=1}^n R_k F_k\right)=\sum_{j=1}^n R_jF_j-\sum_{k=1}^n R_k F_k=0.
\end{eqnarray*}
This says that $\mathcal{P}_xF(x)$ is divergence-free.  Also, note that $\mathcal{P}_x\mathcal{P}_x=\mathcal{P}_x$.  It is enough to prove this on a fixed component.  Again, this follows from direct computation since
\begin{eqnarray*}
\mathcal{P}_x\mathcal{P}_x F_j & = & F_j+R_j\left(\sum_{k=1}^n R_kF_k\right)+R_j\left(\sum_{k=1}^nR_k\left(F_k+R_k\left(\sum_{l=1}^n R_lF_l\right)\right)\right)\\
 & = & F_j+R_j\left(\sum_{k=1}^n R_kF_k\right)+R_j\left(\sum_{k=1}^n R_kF_k\right)+R_j\left(\sum_{k=1}^n R_k^2\left(\sum_{l=1}^n R_lF_l\right)\right)\\
 & = & F_j+R_j\left(\sum_{k=1}^n R_kF_k\right)+R_j\left(\sum_{k=1}^n R_kF_k\right)-R_j\left(\sum_{k=1}^n R_kF_k\right)\\
 & = & F_j+R_j\left(\sum_{k=1}^n R_kF_k\right)=\mathcal{P}_xF_j.
\end{eqnarray*}

Next, we claim that if $b\in BMO(\R^n)$ that we have $[b,\mathcal{P}_x]$ is a bounded operator on $L^2(\R^n;\R^n)$.  Indeed, applying $[b,\mathcal{P}_x]$ to the vector-field $F$ one can easily see that the $j^{th}$ component is given by
$$
\sum_{k=1}^n[b,R_jR_k](F_k)(x).
$$
Indeed, since we have
\begin{eqnarray*}
[b,\mathcal{P}_x](F)_j & = & bF_j+bR_j\left(\sum_{k=1}^n R_kF_k\right)-\left(bF_j+R_j\left(\sum_{k=1}^nR_k\left(bF_k\right)\right)\right)\\
& = & \sum_{k=1}^nR_jR_k(bF_k)-bR_jR_kF_k\\
& = & \sum_{k=1}^n[b, R_jR_k](F_k).
\end{eqnarray*}

With this observation, then an application of the Coifman, Rochberg and Weiss theorem demonstrates $[b,\mathcal{P}_x]$ is bounded on $L^2(\R^n;\R^n)$ with norm controlled by $\norm{b}_{BMO(\R^n)}$.  Then the Div-Curl Lemma follows since one can see using the properties of $\mathcal{P}_x$, and the vector-fields $E(x), B(x)$, that 
\begin{eqnarray*}
\left\vert \int_{\R^n} E(x)\cdot B(x) b(x)dx\right\vert & = & \left\vert \int_{\R^n}[b,\mathcal{P}_x]E(x)\cdot B(x) dx\right\vert\\
& \lesssim & \norm{b}_{BMO(\R^n)}\norm{E}_{L^p(\R^n;\R^n)}\norm{B}_{L^q(\R^n;\R^n)}.
\end{eqnarray*}
But, we already know that this estimate forces $E\cdot B\in H^1(\R^n)$.

We can now use these calculations to handle the proof of the Theorem \ref{uniformly}.  First, observe that we have now two separate projections $\mathcal{P}_x$ and $\mathcal{P}_y$ with each one projecting onto the divergence-free vector fields in the $x$ and $y$ coordinate respectively.  Next, if we fix $y\in\R^n$ then the above computations show that
\begin{eqnarray*}
\int_{\R^n} E(x,y)\cdot B(x,y) b(x,y)dx & = & \int_{\R^n}[b(\cdot,y),\mathcal{P}_x]E(x,y)\cdot B(x,y) dx\\
\end{eqnarray*}
We then integrate this equality over $\R^n$ to see that
\begin{eqnarray*}
\int_{\R^n\times\R^n} E(x,y)\cdot B(x,y) b(x,y)dxdy & = & \int_{\R^n\times\R^n}[b(\cdot,y),\mathcal{P}_x]E(x,y)\cdot B(x,y) dxdy\\
\end{eqnarray*}
 Suppose now that $b\in bmo(\R^n\times\R^n)$.  It is known that a function in this space if for fixed $y$, the function $b(\cdot,y)$ are uniformly in $BMO(\R^n)$.  Thus, the above equality implies
\begin{eqnarray*}
\left\vert \int_{\R^n\times\R^n} E(x,y)\cdot B(x,y) b(x,y)dxdy\right\vert & \leq & \int_{\R^n}\left\vert \int_{\R^n\R^n}[b(\cdot,y),\mathcal{P}_x]E(x,y)\cdot B(x,y) dx\right\vert dy\\
& \lesssim & \int_{\R^n} \norm{b(\cdot, y)}_{BMO(\R^n)} \norm{E(\cdot,y)}_{L^p(\R^n;\R^n)}\norm{B(\cdot,y)}_{L^q(\R^n;\R^n)}dy\\
& \lesssim & \norm{b}_{bmo(\R^n\times\R^n)}\int_{\R^n}\norm{E(\cdot,y)}_{L^p(\R^n;\R^n)}\norm{B(\cdot,y)}_{L^q(\R^n;\R^n)}dy\\
& \lesssim & \norm{b}_{bmo(\R^n\times\R^n)}\norm{E}_{L^p(\R^n\times\R^n;\R^n)}\norm{B}_{L^q(\R^n\times\R^n;\R^n)}.
\end{eqnarray*}
Here we have used the fact that $b\in bmo(\R^n\times\R^n)$ and H\"older's inequality.  Taking the supremum over all $b\in bmo(\R^n\times\R^n)$ we see that
\begin{eqnarray*}
\int_{\R^n}\norm{E(\cdot,y)\cdot B(\cdot,y)}_{H^1(\R^n)}dy\lesssim\norm{E}_{L^p(\R^n\times\R^n;\R^n)}\norm{B}_{L^q(\R^n\times\R^n;\R^n)}.
\end{eqnarray*}
By interchanging the roles of $x$ and $y$ in $\R^n$ we have a similar statement.   This proves Theorem \ref{uniformly}.

\section{Extensions of Theorems \ref{product} and \ref{uniformly} to Differential Forms}

In this section we extend the results of Theorems \ref{product} and \ref{uniformly} to differential forms.  First, we recall some notation.

For $0\leq l\leq n$ the space of $l$-linear, alternating functions $\xi:(\R^n)^{l}\to\R$ is denoted by $\Lambda_{l}(\R^n)$.  The dimension of the space $\Lambda_l(\R^n)$ is $\frac{n!}{l!(n-l)!}$.  The basis for the space $\Lambda_l(\R^n)$ is given by $\{dx_I\}$ where $I$ is an increasing multi-index of $l$ elements, i.e., $I=(i_1,\ldots, i_l)$ with $1\leq i_1<\cdots<i_l\leq n$.

Given $\xi\in\Lambda_l(\R^n)$ and $\eta\in\Lambda_k(\R^n)$ we let $\xi\wedge\eta$ denote the $l+k$ form.  This is given by the following
$$
\xi\wedge\eta\left(X_1,\ldots, X_{l+k}\right)=\sum \textnormal{sgn}(\sigma)\xi\left(X_{i_1},\ldots, X_{i_l}\right)\eta\left(X_{j_1},\ldots, X_{j_k}\right)
$$
where the sum is taken over all possible permutations $\sigma$ of $\{1,\ldots, k+l\}$ such that $i_1<\cdots<i_l$ and $j_1<\cdots<j_k$.  We also let $\textnormal{sgn}(\sigma)$ denote the sign of this permutation.  It is easy to see that $\xi\wedge\eta=-\eta\wedge\xi$, and so the wedge product is alternating.

We will also need the exterior derivative $d$ (sometimes also called the Hodge-de Rham operator).  This operator take $\Lambda_l(\R^n)$ to $\Lambda_{l+1}(\R_n)$ and id defined on functions by 
$$
df:=\sum_{j=1}^n \partial_j f dx_j.
$$
One extends this to general $l$-forms by first acting on the functional coefficient and then simplifying via the exterior algebra induced by the wedge product.  This can then be extended further by using linearity.  A key property will be that $d^2=0$, which is simply the equality of mixed partials.  Finally, let $L^p(\R^n;\Lambda_r)$ denote the space of $p$-integrable $r$-forms on $\R^n$.  One can norm this space in an obvious fashion.

In this language the one-parameter version of Theorem \ref{uniformly} takes the following form.  Suppose that $u\in L^p(\R^n;\Lambda_{n-1})$ with $du=0$, and $v\in L^q(\R^n;\Lambda_1)$ with $dv=0$.  Then $u\wedge v\in H^1(\R^n;\Lambda_n)$.  There is then an extension of this statement to $l$-forms.  Namely, this says that if $1<p<\infty$ and $u\in L^p(\R^n;\Lambda_k)$, $v\in L^q(\R^n;\Lambda_{n-k})$ with $du=0$ and $dv=0$, then $u\wedge v\in H^1(\R^n;\Lambda_n)$.  See \cite{CLMS} for the 1-form case, and \cite{HLMZ} for the case of $k$-forms.

We thus have the following immediate extension of Theorem \ref{uniformly}.

\begin{thm}
\label{uniformlyforms}
Let $1<p<\infty$ and $\frac{1}{p}+\frac{1}{q}=1$.  Suppose that $u\in L^p(\R^n\times\R^n; \Lambda_{k})$ and $v\in L^q(\R^n\times\R^n;\Lambda_{n-k})$ with 
\begin{equation*}
d_x u(x,y)=0\quad d_x v(x,y)=0\quad\forall y\in\R^n
\end{equation*}
and
\begin{equation*}
d_y u(x,y)=0\quad d_y v(x,y)=0\quad\forall x\in\R^n.
\end{equation*}
Then
$$
\int_{\R^n}\norm{u(x,\cdot)\wedge v(x,\cdot)}_{H^1(\R^n)}dx\lesssim\norm{u}_{L^p(\R^n\times\R^n;\Lambda_k)}\norm{v}_{L^q(\R^n\times\R^n;\Lambda_{n-k})}
$$
and
$$
\int_{\R^n}\norm{u(\cdot,y)\wedge v(\cdot,y)}_{H^1(\R^n)}dy\lesssim\norm{u}_{L^p(\R^n\times\R^n;\Lambda_k)}\norm{v}_{L^q(\R^n\times\R^n;\Lambda_{n-k})}.
$$
\end{thm}
The proof of this fact would just be an iteration of the form version of the one-parameter proof and is essentially a repeat of what appear in Theorem \ref{uniformly} so we omit the details.

To extend Theorem \ref{product} we have to introduce a little more notation.  We know are interested in alternating linear functions that map copies of $\R^n\times\R^m$ to $\R$.  For integers $0\leq r\leq n$ and $0\leq s\leq m$ we let $\Lambda_{(r,s)}(\R^n\times\R^m)$ denote the linear space with basis given by $\{dx_I dy_J\}$ where $I$ is an increasing multi-index of length $r$ and $J$ is a multi-index of length $s$.  We can then define the wedge product between the elements $dx_Idy_J$ and $dx_{I'}dy_{J'}$ by 
$$
dx_Idy_J\wedge dx_{I'}dy_{J'}:=dx_I\wedge_x dx_{I'} dy_{J}\wedge_y dy_{J'}.
$$

The generalization to differential forms will then follow in an analogous manner leading to the following Theorem.  

\begin{thm}
\label{productforms}
Let $1<p<\infty$ and $\frac{1}{p}+\frac{1}{q}=1$.  Suppose that $E\in L^p(\R^n\times \R^m;\Lambda_{(r,s)}(\R^n\times\R^m))$ and $B\in L^q(\R^n\times \R^m;\Lambda_{(n-r,m-s)}(\R^n\times\R^m))$ are such that 
\begin{equation*}
d_x E_{I,J}(x,y)=0\quad d_x B_{I',J'}(x,y)=0\quad\forall y\in\R^m,\forall J, J'
\end{equation*}
and
\begin{equation*}
d_y E_{I,J}(x,y)=0\quad d_y B_{I', J'}(x,y)=0\quad\forall x\in\R^n,\forall I, I'.
\end{equation*}
Then $E\cdot B$ belongs to product $H^1(\R^n\otimes\R^m)$ with
$$
\norm{E\wedge B}_{H^1(\R^n\otimes\R^m)}\lesssim\norm{E}_{L^p(\R^n\times\R^m;\Lambda_{r,s}(\R^n\times\R^m)}\norm{B}_{L^q(\R^n\times\R^m;\Lambda_{r,s}(\R^n\times\R^m)}.
$$
\end{thm}



\end{document}